# Comparison between different methods of estimating of the relaxation times in the FPU model

Jacopo De Tullio

Centro PRISTEM, Università commerciale L. Bocconi

April 2015

**Abstract**

After a brief review of the Fermi-Pasta-Ulam (FPU) conservative system of N nonlinearly coupled oscillators, this paper addresses two problems: first, comparing two indicators for the equipartition, showing that the results are essentially identical; second, finding a method that allows fast integration to reach the long integration times required in this area. In particular this work proposes a symplectic algorithm based on the Fast Fourier Transform.

## 1  The FPU model

The first use of a computer to simulate a dynamic system goes back to 1954, when Fermi, Pasta and Ulam conducted a series of numerical experiments that were designed to test the ergodic behavior of a chain of harmonic oscillators linked by not linear springs. Precisely, the Fermi-Pasta-Ulam (FPU) model is formed by a chain of N identical masses interacting with adjacent particles through a weakly nonlinear potential. We define $q_j$ the positions of the particles and $p_j$ the corresponding moments, with $j = 1, \dots, N$ and we also impose the conditions of periodicity to the edge:

$$q_0 = q_N.$$

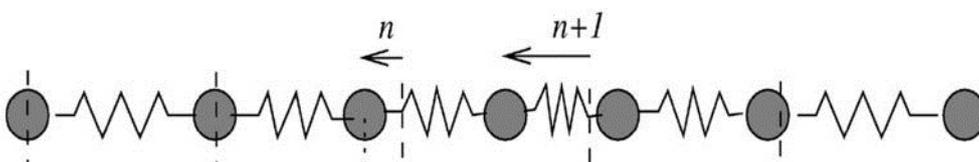

Fig. 1.1: representation of the FPU model; the masses, which can move in only one dimension, interacting nonlinearly in pairs.

The $j - th$ particle the particle will interact only with particles $j - 1$ and $j + 1$, therefore the Hamiltonian of the system just described is given by:

$$H(q_0, q_1, \ldots, p_0, p_1, \ldots, p_N) = \sum_{j=0}^{N-1} \frac{p_j^2}{2} + \sum_{j=0}^{N} V(q_{j+1} - q_j) \quad (1.1)$$

where the mass of the particles was normalized to 1.

The potential $V(q)$ is given by: $V(q) = \frac{1}{2}q^2 + \frac{\alpha}{3}q^3 + \frac{\beta}{4}q^4 \quad \alpha, \beta \geq 0$.

The constants $\alpha, \beta$ are fixed small, in our case they will be $\alpha, = 1/4$ and $\beta = 0$. For $\alpha = \beta = 0$ we obtain a system of equations of linear motion, which can be transformed, by change of variables, in a system of N linear oscillators decoupled. This transformation is defined by:

$$\begin{cases} A_k = \frac{1}{\sqrt{n}} \sum_{j=0}^{N-1} q_j e^{2\pi i k j} \\ \pi_k = \frac{1}{\sqrt{n}} \sum_{j=0}^{N-1} p_j e^{2\pi i k j} \end{cases}$$

that is the discrete transformation of Fourier of $q_j$ and $p_j$.

The transformed Hamiltonian equation is given by:

$$H_2 = \sum_{j=0}^{N-1} \frac{1}{2}(\pi_k^2 + \omega_k^2 A_k^2) \quad (1.2)$$

where the frequency $\omega_k$ are given by $\omega_k = 2 \sin\left(\frac{\pi k}{N}\right)$. The $E_k = \frac{1}{2}(\pi_k^2 + \omega_k^2 A_k^2)$ for $k = 1, \ldots, N$ are constants of motion e are called energy of the normal modes. So we are in the presence of an integrable system because there are N integrals of motion (the harmonic energy $E_k$) independent and in involution (i.e. $\{E_k, E_j\} = 0$ for $k \neq j$).

Considering the anharmonic case, where there is a perturbation in the potential ($\alpha \neq 0$ or $\beta \neq 0$), Fermi, Pasta and Ulam expected that the system were ergodic, i.e. it loses all the first integrals with the exception of the Hamiltonian itself. The passage from harmonic to anharmonic takes to a dichotomy: in the first case there are $N$ independent integrals of motion, in the second one there are no first integrals except the Hamiltonian, so there is the problem of reconciling these two situations with the continuity of the solutions of the equations of the motion with respect to the parameters. A possible solution is given by the notion of relaxation time. It is possible to define as relaxation time the time $\tau$ for which the phase average and the time average essentially coincide.

Returning to the problem of the elimination of the dichotomy between the harmonic and anharmonic case, this is solved by saying that the time $\tau \to +\infty$ for $\alpha, \beta \to 0$. In the case of the FPU is easy to prove that

$$< E_k >_E = \frac{E}{N} \equiv \epsilon \quad k = 1, \ldots, N$$

it means that the expected value of $E_k$, at a given total energy, are all equal independently by $k$. We have in this way the *equipartition of the energy*. The purpose of Fermi, Pasta and Ulam was to define, with a numerical approach, the relaxation time of the time averages $\bar{E}_k(t,x)$ of the energies $E_k$, starting from initial data far from the equilibrium. As initial data they chose the configuration in which the energy is all given on one mode, in particular on the first mode (the one with the lowest frequency):

$$E_1 = E$$
$$E_k = 0 \quad k = 2, \dots, N$$

and with $N = 32$ e the fixed value $\alpha = 1/4, \beta = 1/4, E = 1$.

The expected result was that the energy, due to the nonlinear interactions, would spread over all the other $N - 1$ modes; instead, the results showed that the energy is propagated from the first mode up to the first five and then it returned almost entirely to the first one.
They checked the time averages of the energies $\bar{E}_k(t,x)$ and they noted that the equipartition of energy did not occur but it realized soon a state of apparent equilibrium that was not expected by the hypothesis of ergodicity. These results led to the formulation of the so-called *FPU paradox*: the system instead of having a slow relaxation to the final equilibrium state, it produces a rapid relaxation to a new state in which the energy is divided only between the modes of low frequency, while those of high frequency are excluded. As if the system had an effective number of degrees of freedom lower than $N$.

## 2   Comparison between different methods of estimating

The aim of this section is to estimate the relaxation time, i.e. the time that the FPU system takes to reach the equilibrium. We consider the α-FPU model, i.e. the perturbation is given by the cubic nonlinear interactions; in particular, the Hamiltonian of the system is given by:

$$H(q,p) = \sum_{k=0}^{N} \frac{p_k^2}{2} + \sum_{k=0}^{N+1} \frac{(q_k - q_{k-1})^2}{2} + \alpha \sum_{k=0}^{N+1} \frac{(q_k - q_{k-1})^3}{3} \quad (2.1)$$

with $\alpha = 1/4$.
Will be compared to the times of reaching of equilibrium in two cases. In the first we consider the instantaneous energies of the normal modes of the system; in the second case we build packets of energy of the normal modes, each composed of the same number of components and we will consider the average energy of these packets, defined as the arithmetic average of energies (instantaneous) of the modes that constitute it. In formulas: if $E_k$ is the energy of the $k-th$ normal mode, so $\epsilon_j = \frac{1}{n}\sum_{k=1}^{n} E_{jn+k}$ (in our case $n = 8$) for $j = 0, \dots, \frac{N}{n}$ is the energy of the $j-th$ packet.

### 2.1   Numerical simulation

The numerical simulation of the FPU system is made by a program in language C, written by the author, that implements the algorithm described below. The simulated dynamical system is composed by $N = 512$ particles of mass equal to 1. The input data of the program are:

- initial total energy, which is placed entirely on the normal mode of lower frequency, it assumes various values in the different simulations;

- integration step h for the calculation of the orbits (in the simulations the value set is $h = 2 \cdot 10^{-2}$).

The values of position and time $(q, p)$ of the particles, conjugated to the initial data, are obtained from the normal modes using the inverse Fourier transform. Using the leap-frog method, we proceed to the calculation of the orbit of the system for times of the order of $10^5$.
The total energy of the system is given by:

$$E_{TOT} = \sum_{k=0}^{N} \frac{p_k^2}{2} + \sum_{k=0}^{N+1} \frac{(q_k - q_{k-1})^2}{2} + \alpha \sum_{k=0}^{N+1} \frac{(q_k - q_{k-1})^3}{3} \quad (2.1)$$

we observe that for times of the order of $10^8$, the total energy is saved with an accuracy of $10^{-3}$.

At this point, to calculate the relaxation times, we use an estimator which, in function of time, measures the number of modes that share the energy. At thermodynamic equilibrium the energy will be equally distributed between all modes, so it is expected that the estimator saturates, on an appropriate time scale, to a fixed value that we will calculate. Since the estimator uses the energy of the normal modes, for its calculation it is necessary to pass through *Fast Fourier Transform*, from the variables $(q, p)$ to the normal variables $(A, \pi)$. The instantaneous energy of each modes become:

$$E_j = \frac{\pi_j^2 + \omega_j^2 A_j^2}{2} \quad j = 0, \ldots, N; \quad \omega_j = 2 \sin\left(\frac{\pi j}{N}\right) \quad (2.2)$$

In the simulation the time evolution of the estimator is calculated in two different cases:

- **case 1.** To calculate the estimate we use the instantaneous energy $E_j$ for $j = 0, \ldots, N$.

- **case 2.** Instead of $N$ energies are considered packets of energy $\varepsilon_j$, each consisting of $n$ consecutive modes, with $n$ proportional to $N$: $\varepsilon_j = \sum_{k=nj+1}^{n(j+1)} E_k$, $j = 0, \ldots, N-1$, in particular we have chosen $n = 8$, so $N_0 = 64$.

We study the relaxation time in the two previous case to see if there are relevant difference (as it happens in [4]).

## 2.2 Estimator $n_{eff}$

We start introducing the quantity:

$$S = -\sum_{i=0}^{N} e_i \ln e_i \quad (2.3)$$

where $e_i$ are the normalized instantaneous energy of the normal modes, i.e.

$$e_i = \frac{E_i}{\sum_{i=0}^{N} E_i}$$

Now we consider the normalized estimator $n_{eff} = \frac{e^S}{N}$ (2.4). We can observe that in the case in which there is equipartition of the instantaneous energy $E_i = \frac{E_{TOT}}{N}$ we have $e_i = \frac{1}{N}$, so the entropy becomes:

$$S = -\sum_{i=0}^{N} \frac{1}{N} \ln \frac{1}{N} = -\ln \frac{1}{N} = \ln N$$

and the estimation becomes

$$n_{eff} = \frac{e^{\ln N}}{N} = 1$$

it means that all the modes share the system's energy.

In the simulations the value of $n_{eff}$ doesn't go to 1, it means that the energy doesn't share itself on all the modes because of the energy's fluctuations of each mode. To calculate the effect of the fluctuations we introduce a deviation $\delta e_i$ of the equipartition: $e_i = \bar{e} + \delta e_i$ (2.5) where $\bar{e} = \frac{1}{N}\sum_{i=0}^{N} e_i$ is the average instantaneous energy.

Substituting (2.5) in (2.3) and developing the logarithm with Taylor polynomial as:

$$\ln\left(1 + \frac{\delta e_i}{\bar{e}}\right) = \frac{\delta e_i}{\bar{e}} - \frac{1}{2}\left(\frac{\delta e_i}{\bar{e}}\right)^2$$

we have:

$$S = -\sum_{i=0}^{N} (\bar{e} + \delta e_i) \ln(\bar{e} + \delta e_i) = -\sum_{i=0}^{N} (\bar{e} + \delta e_i) \ln\left[\bar{e}\left(1 + \frac{\delta e_i}{\delta e_i}\right)\right] =$$

$$= -\sum_{i=0}^{N} (\bar{e} + \delta e_i) \ln \bar{e} - \sum_{i=0}^{N} (\bar{e} + \delta e_i)\left[\frac{\delta e_i}{\bar{e}} - \frac{1}{2}\left(\frac{\delta e_i}{\bar{e}}\right)^2\right] =$$

$$= -\sum_{i=0}^{N} \bar{e} \ln \bar{e} - \ln \bar{e} \sum_{i=0}^{N} \delta e_i - \sum_{i=0}^{N} (\bar{e} + \delta e_i)\frac{\delta e_i}{\bar{e}}\left(1 - \frac{\delta e_i}{2\bar{e}}\right) =$$

if we use the fact that $\sum_{i=0}^{N} \delta e_i = 0$ because of the sum of all the deviations is null, we have:

$$S = -N\bar{e}\ln\bar{e} - \sum_{i=0}^{N}\left(\delta e_i - \frac{(\delta e_i)^2}{\bar{e}}\right)\cdot\left(1 - \frac{\delta e_i}{2\bar{e}}\right) = -N\bar{e}\ln\bar{e} - \sum_{i=0}^{N}\left(\delta e_i + \frac{(\delta e_i)^2}{2\bar{e}}\right) =$$

$$= -N\bar{e}\ln\bar{e} - \frac{1}{2\bar{e}}\sum_{i=0}^{N}(\delta e_i)^2$$

At the equilibrium the Boltzmann's is valid and it implies $(\delta e_i)^2 = (\delta \bar{e})^2$. So we have:

$$-N\bar{e}\ln\bar{e} - N\frac{(\delta \bar{e})^2}{2\bar{e}} \quad (2.6)$$

Substituting (2.6) in (2.4) we have:

$$n_{eff} = \frac{1}{N}\exp\left\{-N\bar{e}\ln\bar{e} - N\frac{(\delta\bar{e})^2}{2\bar{e}}\right\} = \frac{1}{N}\exp(-N\bar{e}\ln\bar{e}) \cdot \exp\left(-N\frac{(\delta\bar{e})^2}{2\bar{e}}\right) =$$

$$= \exp\left(-N\frac{(\delta\bar{e})^2}{2\bar{e}}\right) \quad (2.7)$$

Because of $\bar{e} = \frac{1}{N}$ and $(\delta\bar{e})^2 = \bar{e}^2$, we obtain the asymptotic value of the estimation: $n_{eff} = \exp\left(-\frac{1}{2}\right) = 0{,}61$.

These calculations show that:

- the result doesn't depend on the number of the oscillators if $N$ is sufficiently high;
- the limit value of the estimator for $n \to +\infty$.

## 2.3 Numerical results

Proceeding with the numerical simulations we can determine the relaxation times of the system and compare the values provided by the two estimators. Recalling that the chain is composed by 512 particles and all the energy is placed on the first mode.

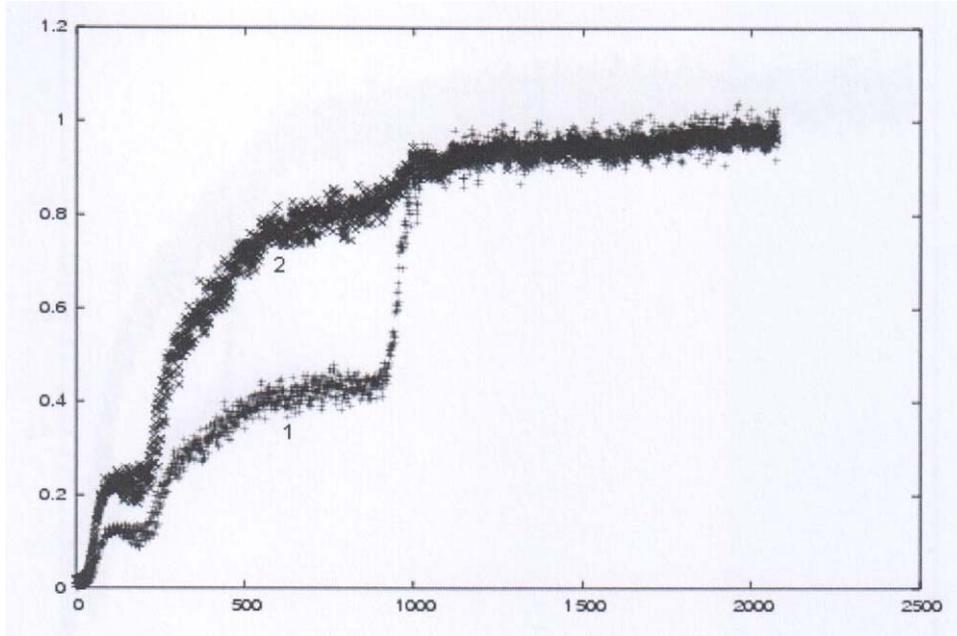

Fig. 2.1: comparison of the two estimates $n_{eff}$ for the calculation of the times of equilibrium. The first estimate is made with the instantaneous energies, the second with the packets of energy. Input value $A_1 = 40$, in particular the values are normalized to 1.

The Fig. 2.1 shows the comparison between the two estimators of relaxation times made on the FPU model in which we have assigned all the energy to the first mode (in particular $A_1 = 40$). We observe that initially the estimate made with the packets seems to be faster, while at the end both reach the equipartition at time $T_{eq} \approx 2 \cdot 10^3$. About the time $10^3$ the trend of the estimate calculated with the instantaneous energies makes a jump coming to the same values of the estimate on the packets and the two curves almost overlap.

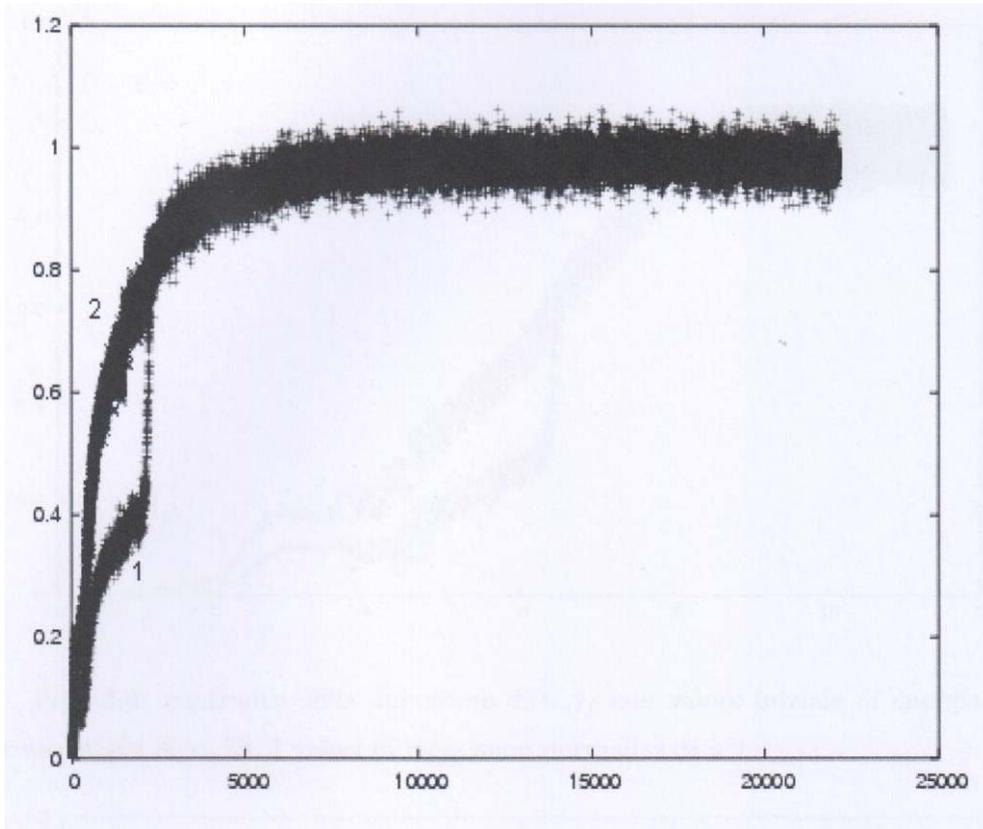

Fig. 2.2: input value $A_1 = 30$.

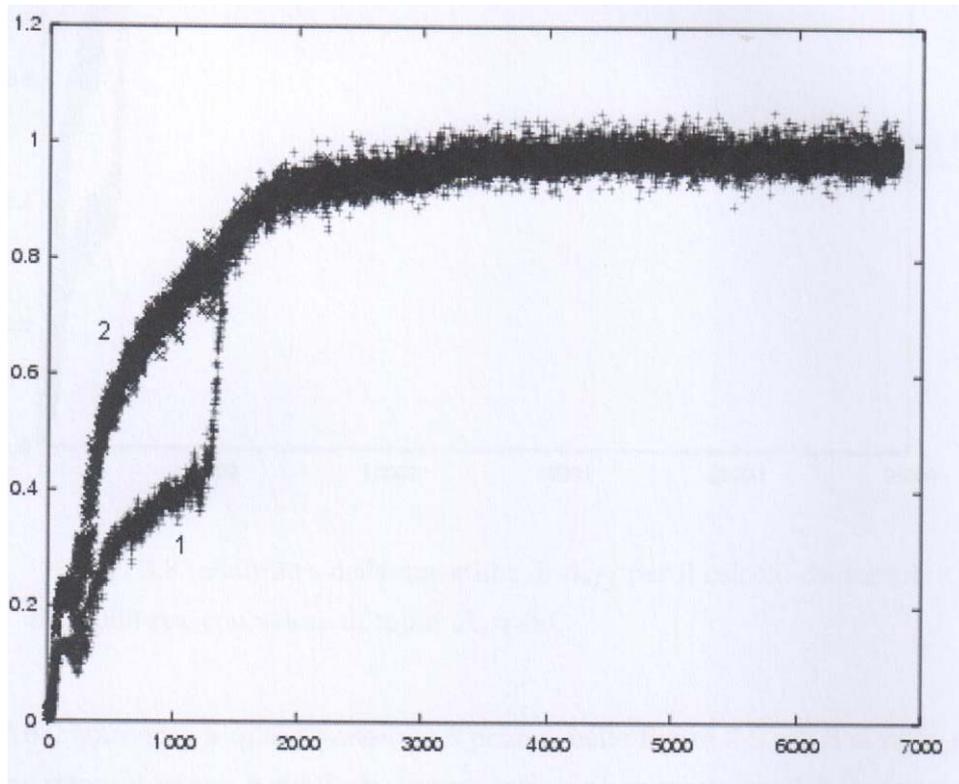

Fig. 2.3: input value $A_1 = 35$.

Similarly to what was observed before, in Fig. 2.2 and Fig. 2.3 we see that the two estimates lead to equivalent results. In particular in the case $A_1 = 35$ equilibrium times are of the order of $3 \cdot 10^3$, while in the case $A_1 = 30$ times are of the order $6 \cdot 10^3$.

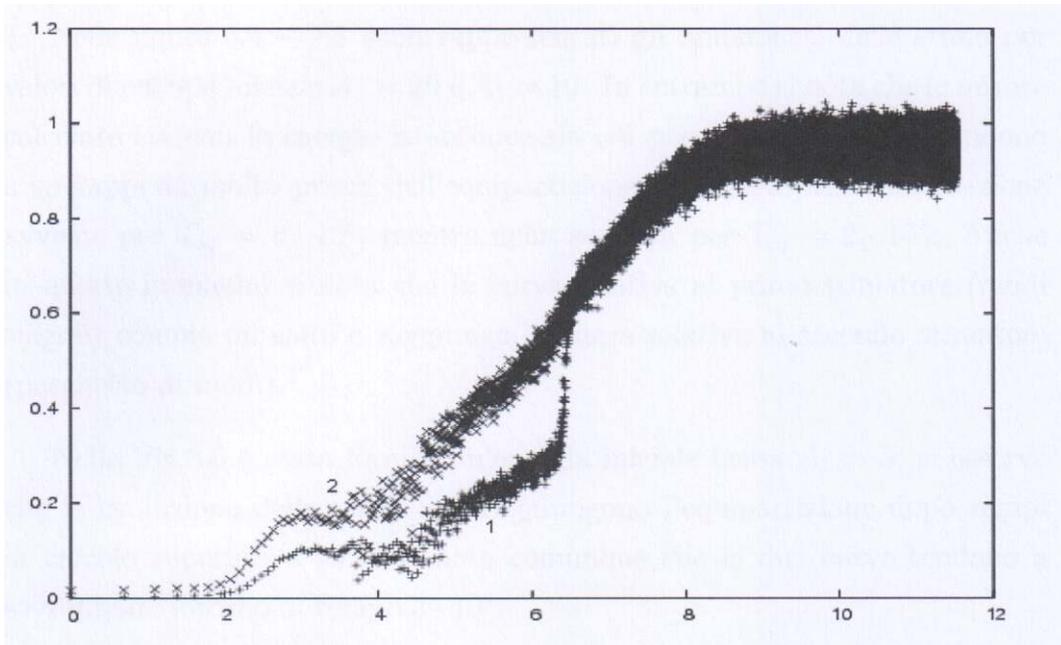

Fig. 2.4: input value $A_1 = 20$.

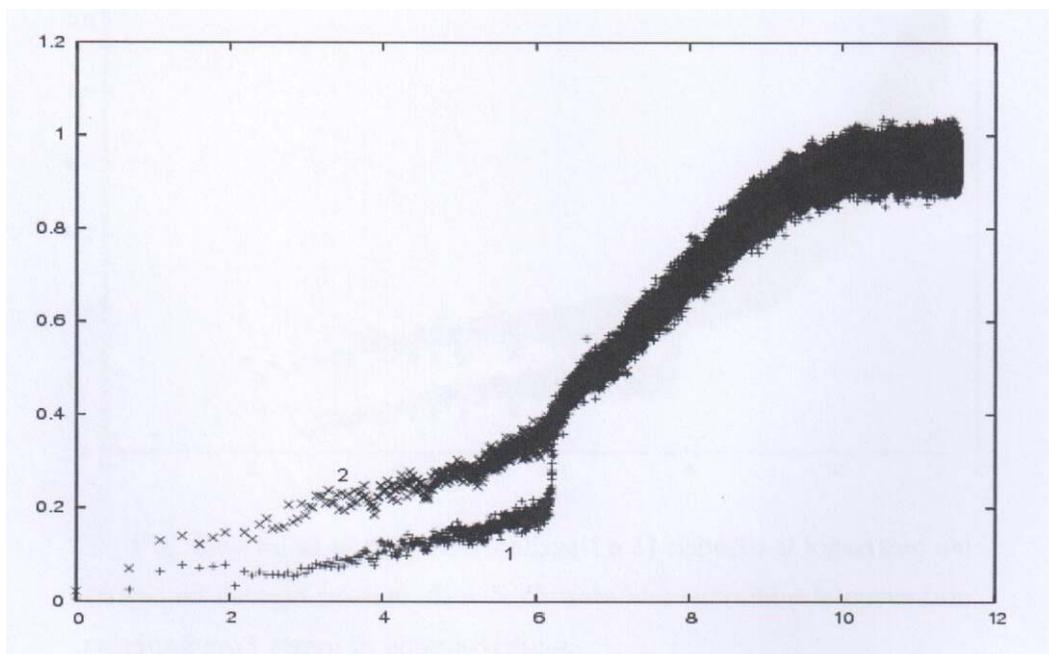

Fig. 2.5: input value $A_1 = 10$.

In Fig. 2.4 and Fig. 2.5 are shown the estimates' trends for values of initial energy $A_1 = 20$ and $A_1 = 10$. In both the case the curves overlap long before the equipartition, that is obtaned at $T_{eq} \approx 8 \cdot 10^3$ (first case) and at $T_{eq} \approx 2 \cdot 10^4$ (second case). In Fig. 2.6 the initial energy given is very low ($A_1 = 5$); we observe that the evolution of the estimates don't reach the equipartition for calculation times higher than $10^5$.

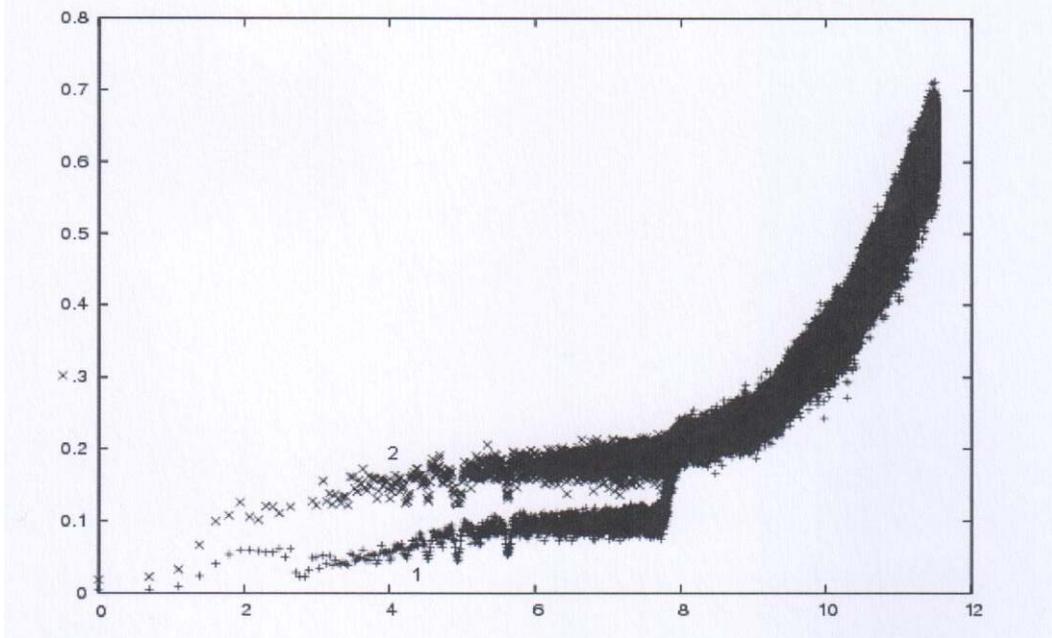

Fig. 2.6: input value $A_1 = 5$. The value of $n_{eff}$ are normalized with respect to the logarithm of time.

From the data it is deduced that the two estimates considered and implemented give results almost equal to the calculation of the system's times of equilibrium. We note that the performance of the estimate computed through packets of energy turns out to have less variation in the values that it takes and is more regular than the estimate given by the instantaneous energy that has jumps and changes of convexity. In addition, from the values of the times of equilibrium we observe that to lower energies correspond long time to achieve the equipartition.

# 3 New method of integration

This section proposes a fast method of integration that allow to reach the long times required for the equipartition.
The method used in the preview section (the leap-frog method) shows that if we start with low energy it reaches long times to reach the equipartition, also using integration time of the order of $10^8$ the system stays far from the equilibrium. For this reason we built an alternative method of integration, fast and efficient, to calculate the orbits of the system also for very long times. This new algorithm, adopted for the α-FPU, model is based on a linearized system.

## 3.1 Alternative method of integration

Consider the FPU's Hamiltonian

$$H(q,p) = \underbrace{H_0}_{quadratic} + \underbrace{H_1}_{cubic} = T + V$$

where $V = V_2 + V_3$ is the sum among the quadratic and cubic potential. In the leap-frog method $H_0 = T$ and $H_1 = V$. Now we give a different way considering $H_0(p,q) = T + V_2$ and $H_1(q) = V_3$. We note that $H_1$ depends only by $q$ so it's easy to find its flux $\Phi_1^h$:

$$\begin{cases} q^{(n+1)} = q^{(n)} \\ p^{(n+1)} = p^{(n)} - h\alpha \left[ \left(q_k^{(n)} - q_{k-1}^{(n)}\right)^2 - \left(q_{k+1}^{(n)} - q_k^{(n)}\right)^2 \right] \end{cases} \quad (3.1)$$

Then we convert $H_0(p,q)$ into normal modes trough Fourier transform:

$$H_0(A, \pi) = \sum_{k=0}^{N} \frac{\pi_k^2 + \omega_k^2 A_k^2}{2}, \quad \omega_k = 2\sin\left(\frac{\pi k}{N}\right)$$

At this point the flux $\Phi_0^h$ of the Hamiltonian $H_0$ is given by:

$$\begin{cases} \pi_k^{(n+1)} = \pi_k^{(n)} \cos(\omega_k \cdot h) - \omega_k A_k^{(n)} \sin(\omega_k \cdot h) \\ A_k^{(n+1)} = A_k^{(n)} \cos(\omega_k \cdot h) - \frac{\pi_k^{(n)}}{\omega_k} \sin(\omega_k \cdot h) \end{cases} \quad (3.2)$$

where the flux is given by the solutions of the harmonic oscillators. In fact it's easy to prove that (3.2) satisfies the Hamilton equation

$$\begin{cases} \dot{A}_k = \frac{\partial H_0}{\partial \pi_k} = \pi_k \\ \dot{\pi}_k = -\frac{\partial H_0}{\partial A_k} = -\omega_k^2 A_k \end{cases}$$

so, through the flux $\Phi_0^h$, we can compute the orbits of the linear part of the system. Now it is easy to check that the map $\Phi_1^{h/2} \circ \Phi_0^h \circ \Phi_1^{h/2}$ approximates the real flux $\Phi^h$ of the system with an error of order $h^3$.

## 3.2 Numerical simulations

The numerical simulation of the system was made by a program in C language (made by the author) that as input has the values of the normal modes and the step of integration. Similarly to the previous case we give all the energy to the lowest frequency mode. We compute the orbits of the linear part using the alternative method just exposed and for the nonlinear part we use the leap-frog.
The estimators $n_{eff}$ are calculated in two cases: by the instantaneous energy of the normal modes and by the packets of normal modes. In particular we have simulated a FPU system made by *N=512* particles with step of integration $h = 1$ and a low initial energy ($A_1$=2 or 5). So, using a step of integration 100 times bigger than in the previous simulations, we have reached times of integration of order $t = 10^{10}$. Despite we have increased the integration step, we preserve the total energy of the system with an accuracy of $10^{-4}$. We also observe that the new algorithm is faster than the previous, so it is possible to increase the number of the orbits that we want calculate without that the calculators uses too much time.

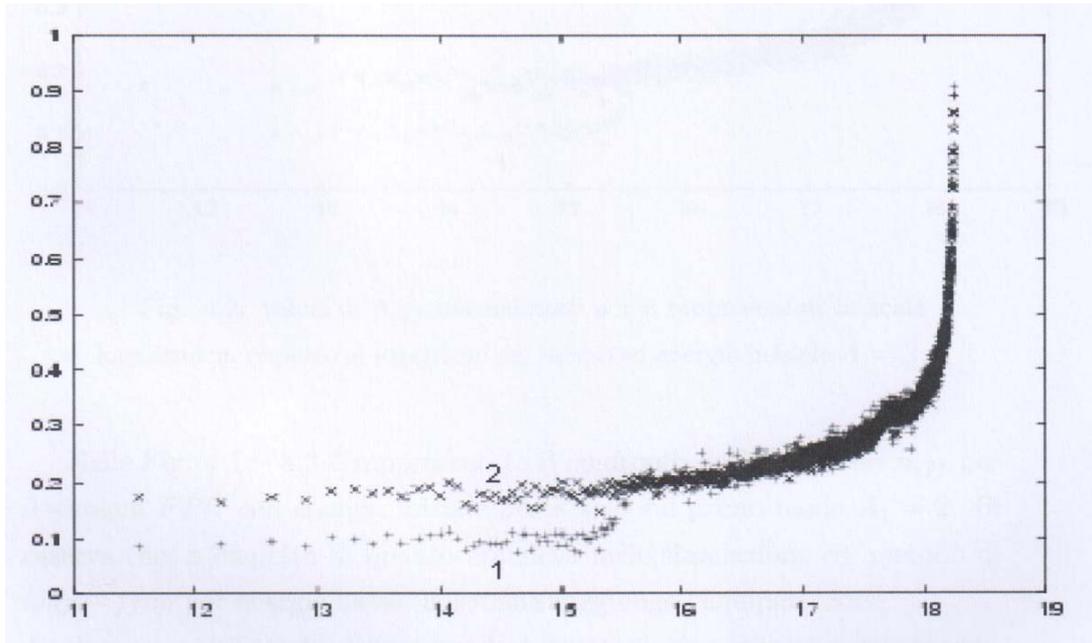

Fig. 3.1: comparison of the two estimates $n_{eff}$ for the calculation of the equilibrium times (in logarithmic scale). The first estimate is made with the energies snapshots with the 2 packets of energy. Input values $A_1$=2 and time $T = 10^{10}$.

The figures 3.1 shows the comparison of the estimator $n_{eff}$ where the initial energy in on the first mode $A_1 = 2$. We observe that, despite what happened with the leap-frog method, the system reach the equipartition. We also note that the two estimators have the same behavior and give similar results for the determination of the equilibrium time.

This alternative method provides results similar to leap-frog, with the advantage of using much smaller computing time and preserves more of the total energy. In this way we have an algorithm that allow to calculate accurately and with greater integration times the equilibrium time of the system FPU in the case of low initial energy.

# References

[1] Benettin G., Carati A., Galgani L., Giorgilli A.: The Fermi-Pasta-Ulam problem and the metastability perspective, in "The Fermi-Pasta-Ulam Problem: A Status Report", G. Gallavotti editor, Lecture Notes in Physics , Vol. 728, Springer Verlag (Berlin,2007).

[2] Fermi E., Pasta J., Ulam S.: Studies of Non Linear Problems, Document LA-1940 (1955).

[3] De Luca J., Lichtenberg A.: Transitions and time scale to equipartition in oscillator chaines: Low-frequency initial conditions, *Phys. Rev. E*, 66, 026206 (2002).

[4] Carati A., Galgani L., Giorgilli A., Paleari S.: FPU phenomenon for generic initial data, *Phys. Rev. E*, 76, 022104 (2007).

[5] Maiocchi A., Bambusi D., Carati A.: An Averaging Theorem for FPU in the Thermodynamic Limit, *J. Stat. Phys.*, 155, 300 (2014).